\documentclass[10pt,a4paper]{amsart}  
\usepackage{a4wide,amsmath,amssymb,mathptmx,hyperref}
\usepackage[english,francais]{babel}

\newtheorem{theorem}{Theorem}[section]
\newtheorem{lemma}[theorem]{Lemma}
\newtheorem{e-proposition}[theorem]{Proposition}

\newtheorem{e-definition}[theorem]{Definition\rm}
\newtheorem{remark}{\it Remark\/}

\providecommand{\BBb}[1]{{\mathbb{#1}}}
\providecommand{\cal}[1]{{\mathcal{#1}}}   

\newcommand{\B}{{\BBb B}}
\newcommand{\C}{{\BBb C}}
\newcommand{\im}{\operatorname{i}}
\newcommand{\nrm}[2]{\|#1\|_{#2}}
\newcommand{\Nrm}[2]{\bigl\|#1\bigr\|_{#2}}
\newcommand{\norm}[2]{\mathinner{\|}#1\,|#2\|}
\newcommand{\Norm}[2]{\mathinner{\bigl\|\,#1\,\big|#2\bigr\|}}
\newcommand{\N}{\BBb N}
\newcommand{\op}[1]{\operatorname{#1}}
\newcommand{\R}{{\BBb R}}
\newcommand{\Rn}{{\BBb R}^{n}}
\newcommand{\supp}{\operatorname{supp}}
\renewcommand{\hat}[1]{\overset{{\scriptscriptstyle \wedge}}{#1}}

\begin{document} 
\title[Operators of type $1,1$ and Triebel--Lizorkin spaces]{Domains of type $1,1$ operators: a case for 
Triebel--Lizorkin spaces
\\
Domaines des op\'erateurs de type $1,1$ et espaces de
Triebel--Lizorkin}
\author{Jon Johnsen}
\address{Department of Mathematics, Aalborg University,
Fredrik Bajers Vej 7G; DK-9220 Aalborg {\O}st, Denmark.}
\email{\tt jjohnsen@math.aau.dk}
\thanks{~
\\[4\jot] {\tt Appeared in Comptes Rendus Academie de Sciences Paris, Serie I, {\bf 339} (2004), no.~2,
   115--118.}}
\maketitle

{\small
{\bf Abstract.\/}
Pseudo-differential operators of
type $1,1$ are proved continuous from the Triebel--Lizorkin space
$F^d_{p,1}$ to $L_p$, $1\le p<\infty$, when of order $d$, 
and this is in general 
the largest possible domain among the Besov and Triebel--Lizorkin
spaces. H{\"o}rmander's condition on the twisted diagonal is extended to
this framework, using a general support rule for Fourier
transformed pseudo-differential operators.%
}

{\small 
{\bf R\'esum\'e.\/}
On d\'emontre que les op\'erateurs pseudo-diff\'erentiels de type $1,1$ et 
d'ordres $d$ sont
continus de l'espace $F^d_{p,1}$ de Triebel--Lizorkin dans $L_p$, $1\le
p<\infty$, et que parmi les espaces de Besov et Triebel--Lizorkin,
ces domaines sont en g\'en\'eral le plus grand possible.
La condition de
H{\"o}rmander sur la diagonale tordu est \'etablie pour ce cadre, 
en utilisant un
resultat g\'en\'eral sur le support de la transformation de Fourier d'un
operateur pseudo-diff\'erentiel.
}
\section{Introduction}   \label{intr-sect}
\enlargethispage{\baselineskip}
\thispagestyle{empty}

Recall that for symbols $a\in S^{d}_{\rho,\delta}(\Rn\times\Rn)$, ie
$|D^\alpha_\xi D^\beta_x a(x,\xi)|\le 
   C_{\alpha\beta} (1+|\xi|)^{d-\rho|\alpha|+\delta|\beta|}$,
\begin{equation}
  a(x,D)=\op{OP}(a)=(2\pi)^{-n}\!\int 
  e^{\im x\cdot\xi} a(x,\xi)\hat u(\xi)\,d\xi
  \label{axD-eq}
\end{equation}
map the Schwartz space $\cal S(\Rn)$ continuously into itself,
say for $0\le \delta\le\rho\le1$. 
And for $(\rho,\delta)\ne (1,1)$ these operators extend to continuous,
`globally' defined maps
\begin{equation}
  a(x,D)\colon \cal S'(\Rn)\to\cal S'(\Rn).
  \label{aS-eq}
\end{equation}
But for $\rho=\delta=1$ Ching \cite{Chi72} proved existence of
$a\in S^0_{1,1}$ such that $a(x,D)\notin\B(L_2(\Rn))$. 
That every $A\in\op{OP}(S^0_{1,1})$ is bounded
on $C^s$ and $H^s$ for $s>0$ was first proved by Stein
(unpublished);
Meyer \cite{Mey80} proved 
continuity from $H^{s+d}_p$ to $H^{s}_p$ for $s>0$, $1<p<\infty$.

For $s\le0$, H{\"o}rmander \cite{H88} gave a condition 
on the \emph{twisted diagonal} $\{\,(\xi,\eta)\mid \xi+\eta=0\,\}$: 
$a(x,D)$ is bounded $H^{s+d}\to H^{s}$ for \emph{all} $s\in\R$
if $\hat a(\xi,\eta):=\cal F_{x\to\xi}a(x,\eta)$ fulfils 
\begin{equation}
  \hat a(\xi,\eta)=0 \quad\text{for}\quad 
  C(|\xi+\eta|+1)\le |\eta|,\quad\text{for some}\quad C\ge1.
  \label{H-cnd}
\end{equation}

For $s\ge0$ and $1\le p\le\infty$, the next result gives a maximal domain  
by means of the Triebel--Lizorkin spaces $F^{s}_{p,q}(\Rn)$ (albeit with a
Besov space for $p=\infty$). 

\begin{theorem}
  \label{Fp1-thm}
Every $a\in S^{d}_{1,1}(\Rn\times\Rn)$, $d\in\R$, 
gives a bounded operator 
\begin{align}
  a(x,D)&\colon F^{d}_{p,1}(\Rn)\to L_p(\Rn)
         \quad\text{for}\quad p\in[1,\infty[,
\\
  a(x,D)&\colon B^{d}_{\infty,1}(\Rn)\to L_\infty(\Rn).
  \label{Fp1-eq}
\end{align}
$\op{OP}(S^d_{1,1})$ contains  $a(x,D)\colon \cal S(\Rn)\to\cal
D'(\Rn)$, that are discontinuous when $\cal S(\Rn)$ is given the induced
topology from any $F^d_{p,q}(\Rn)$
or $B^d_{p,q}(\Rn)$ with $p\in [1,\infty]$ and $q\in\,]1,\infty]$.
\end{theorem}

So for fixed $p\in[1,\infty[\,$, every
$A\in\op{OP}(S^0_{1,1})$ is bounded 
$F^0_{p,1}\to L_p$ and everywhere defined, 
but not so on any larger $B^{s}_{p,q}$- or 
$F^{s}_{p,q}$-space (regardless of the codomain).

In comparison with Besov spaces,
arguments in favour of \emph{Triebel--Lizorkin} spaces have,
perhaps, been less convincing. Indeed, $F^s_{p,2}=H^s_p$ for
$1<p<\infty$, cf.~\cite{T2}, but this doesn't necessarily make
the $F^s_{p,q}$ a \emph{useful} extension of 
the $H^s_p$-scale.
However, Theorem~\ref{Fp1-thm} shows that also $F^{s}_{p,q}$-spaces with
$q\ne 2$  are indispensable for a natural $L_p$-theory. 

The next result extends
H{\"o}rmander's condition in \eqref{H-cnd} to $F^{s}_{p,q}$.

\begin{theorem}
  \label{cont-thm}
Any $a(x,D)\in \op{OP}(S^{d}_{1,1}(\Rn\times\Rn))$ is continuous, 
for $s>0$, $p$, $q\in[1,\infty]$,
\begin{equation}
  a(x,D)\colon F^{s+d}_{p,q}(\Rn)\to F^{s}_{p,q}(\Rn), 
    \quad\text{for}\quad p<\infty.
  \label{Fs-eq}
\end{equation}
If \eqref{H-cnd} holds, \eqref{Fs-eq} does so for $s\in\R$.
(The result extends to $B^{s}_{p,q}$ and $p,q\in\,]0,\infty]$).
\end{theorem}

The proofs of Theorem~\ref{Fp1-thm}--\ref{cont-thm} treat
the symbols directly without approximation by elementary symbols, so
it is crucial to control the spectra of the terms appearing in the
paradifferential splitting of $a(x,D)$, 
and for this purpose the following was established.

\begin{e-proposition}[the support rule]
  \label{supp-prop}
If $b\in S^d_{1,0}(\Rn\times\Rn)$ and $v\in\cal F^{-1}\cal
E'(\Rn)$, then
\begin{equation}
  \supp \cal F(b(x,D)v) \subset
  \bigl\{\,\xi+\eta \bigm| (\xi,\eta)\in \supp \hat b(\cdot,\cdot),
            \eta\in\supp \hat v \,\bigr\}.
  \label{supp-eq}
\end{equation}
\end{e-proposition}

\begin{e-proposition}   \label{E'-prop}
  Any $A$ in $\op{OP}(S^\infty_{1,1})$ extends to a map $\cal F^{-1}\cal
E'(\Rn)\to \cal S'(\Rn)$, that coincides with the usual one for 
$A\in \op{OP}(S^\infty_{1,0})$.
\end{e-proposition}

The support rule generalises to $b\in S^\infty_{1,1}$, for all $v\in\cal
F^{-1}\cal E'$, using Proposition~\ref{E'-prop}.

\section{On the proofs}
  \label{brdl-sect}
With $1=\sum_{j=0}^\infty\Phi_j$ so that $\Phi_j(\xi)=1\iff 
|\xi|\sim 2^j$ ($j>0$), set 
$\tilde\Phi_j=\Phi_{j-1}+\Phi_j+\Phi_{j+1}$,
$a_{j,k}(x,\eta)=\cal F^{-1}_{\xi\to x}(\Phi_j\hat a(\cdot,\eta))
\tilde\Phi_k(\eta)$ and $u_j=\Phi_j(D)u$. 
One can then make the ansatz
\begin{equation}
  a(x,D)u(x)=a^{(1)}(x,D)u(x)+a^{(2)}(x,D)u(x)+a^{(3)}(x,D)u(x),
  \label{a123-eq}
\end{equation}
when the pair $(a,u)$ is such that the following
series converge in $\cal D'(\Rn)$:
\begin{align}
  a^{(1)}(x,D)u&=\sum_{k=2}^\infty \sum_{j=0}^{k-2} a_{j,k}(x,D)u_k, 
\qquad  a^{(3)}(x,D)u=\sum_{j=2}^\infty \sum_{k=0}^{j-2} a_{j,k}(x,D)u_k \\
  a^{(2)}(x,D)u&=\sum_{k=0}^\infty \sum_{j,l=0,1,\ j+l\le1}
                  a_{k-j,k-l}(x,D)u_{k-l}, 
\end{align}
Here $a\in S^\infty_{1,1}(\Rn\times\Rn)$ implies $a_{j,k}\in S^{-\infty}$,
and if $K_{j,k}$ denotes the distribution kernel,
\begin{equation}
  a_{j,k}(x,D)u_k=\int_{\Rn}K_{j,k}(x,y)u_k(y)\,dy, 
   \quad\text{for}\quad u\in\cal S'(\Rn).
   \label{ajk-eq}
\end{equation}
This definition of $a(x,D)$ extends other ones, eg \eqref{axD-eq}.
And Prop.~\ref{E'-prop} follows, for
if $\hat u\in\cal E'$ both $a^{(1)}(x,D)u$, $a^{(2)}(x,D)u$ exist
as finite sums; with 
$K_k(x,y):=\cal F^{-1}_{\xi\to y}(a\tilde\Phi_k)(x,x-y)$
one can sum over $j\le N$ in \eqref{ajk-eq} and majorise to show 
$\cal S'$-convergence to $\int K_k(x,\cdot)u_k\,dy$.

To exploit the ansatz further, the `pointwise' estimate in the next lemma is
useful.
\begin{lemma}
  \label{JM-lem}
Let $v\in \cal S'(\Rn)$ and $b\in
S^\infty_{1,1}(\Rn\times\Rn)$ such that 
$\supp\cal Fv\cup\bigcup_{x\in\Rn}\supp b(x,\cdot)$ is contained in a
ball $B(0,2^k)$, $k\in\N$.
Then there exists a $c>0$ such that 
\begin{equation}
  |b(x,D)v(x)|\le c\norm{b(x,2^k\cdot)}{\dot B^{n/t}_{1,t}(\Rn)}
  M_tv(x).
  \label{JM-ineq}
\end{equation} 
Here 
$M_tf(x)=\sup_{r>0}(\tfrac{1}{|B(x,r)|}\int_{B(x,r)}|f(y)|^t\,dy)^\frac1t$ 
is the maximal function; $0<t\le 1$.
\end{lemma}

Lemma~\ref{JM-lem} is similar to
\cite[Prop.~5(a)]{Mar96}, except that
$b\in S^\infty_{1,1}$ replaces the vague
assumption of being a `symbol $\Rn\times\Rn\to\C$'
(\cite[Prop.~5(a)]{Mar96} itself is not easy to read, 
as it is extracted from an
earlier proof with another set-up. But $b\in S^\infty_{1,1}$ 
implies that $b(x,D)v$ is given by an integral like \eqref{ajk-eq},
and estimates in \cite[Prop.~4]{Mar96} apply to this.)

The proof of Theorem~\ref{Fp1-thm} combines \eqref{JM-ineq} with 
$L_p(\ell_1)$-boundedness of $M_t$ for $t<1$, so that $\tfrac{n}{t}<n+1$. 
Further estimates of $a$ follow from the embeddings 
$W^{n+1}_1\hookrightarrow B^{n+1}_{1,\infty}\hookrightarrow 
\dot B^{n/t}_{1,t}$: since $\tfrac{1}{4}\le|\eta|\le4$ on $\supp\tilde\Phi$,
so eg $2^{kd}\sim (1+|2^k\eta|)^d$, then if $\Psi_k=\Phi_0+\dots+\Phi_k$,
\begin{equation}
  \Norm{\sum_{j=0}^{k-2}a_{j,k}(x,2^k\cdot)}{\dot B^{n/t}_{1,t}}
 \le 
  \sum_{|\alpha|\le n+1}\Norm{D^\alpha_\xi(\Psi_{k-2}(D_x)a(x,2^k\cdot)
    \tilde\Phi)}{L_{1,\xi}}
 \le c 2^{kd},
\end{equation}
where $c=c'\norm{\tilde\Phi}{W^{n+1}_1}
              \nrm{\check\Psi}{1}\sup_{x,\xi; |\alpha|\le n+1}
 (1+|\xi|)^{-(d-|\alpha|)}|D^\alpha_\xi a(x,\xi)|$.
Using \eqref{JM-ineq},
\begin{equation}
  \begin{split}
    \Nrm{\sum_k\sum_{j=0}^{k-2}a_{j,k}(x,D)u_k}{p}^p
  &\le \int |\sum_k 2^{kd} M_t u_k(x)  |^p\,dx
 (\sup_{x,k}2^{-kd}
  \Norm{\sum_{j=0}^{k-2}a_{j,k}(x,2^k\cdot)}{\dot B^{n/t}_{1,t}})^p
\\
  &\le c \int (\sum_k 2^{kd}|u_k(x)|)^p\,dx.
  \end{split}
  \label{a1-ineq}
\end{equation}
For $k$ in finite sets, it now follows that the $a^{(1)}(x,D)u$-series
is fundamental in $L_p$ when $u\in F^d_{p,1}(\Rn)$ for $1\le p<\infty$, and
\eqref{a1-ineq} gives that $a^{(1)}(x,D)$ is bounded.
The sum $\sum_{j=0}^{k-2}$ may then be replaced by the one
pertinent for $a^{(2)}$, with a similar argument.
To handle $a^{(3)}$, one may further invoke Taylor's formula and
\cite[Lem.~3.8]{Y1}. 
The case $B^d_{\infty,1}(\Rn)$ is analogous, and the 
counterexamples of \cite{Chi72} adapts easily to give the sharpness.

In the proof of Theorem~\ref{cont-thm},
the key point is to obtain (with $\Phi_j$ as in \cite{Y1})
\begin{gather}
  \supp\cal F\bigl(\sum_{j=0}^{k-2} a_{j,k}(x,D)u_k\bigr)
\cup  \supp\cal F\bigl(\sum_{j=0}^{k-2} a_{k,j}(x,D)u_j\bigr)
  \subset
  \bigl\{\, \tfrac{1}{5}2^k\le 
                          |\xi|\le 5\cdot2^k\,\bigr\},
   \label{S1}\\
  \supp\cal F\bigl(\sum_{j,l=0,1\ j+l\le1} 
                         a_{k-j,k-l}(x,D)u_{k-l}\bigr)
  \subset
  \bigl\{\, |\xi|\le 4\cdot2^k\,\bigr\}.
  \label{S2}
\end{gather}
If \eqref{H-cnd} holds,
then \eqref{S2} may be supplemented by the property that, 
for $k$ large enough, 
\begin{equation}
    \supp\cal F\bigl(\sum_{j,l=0,1\ j+l\le1} 
                         a_{k-j,k-l}(x,D)u_{k-l}\bigr)
  \subset
  \bigl\{\,\xi \bigm| \tfrac{1}{4C}2^k\le
                          |\xi|\le 4\cdot2^k\,\bigr\}.
  \label{S2'}
\end{equation}
By Proposition~\ref{supp-prop}, \eqref{S1}--\eqref{S2} are easy.
\eqref{S2'}
is seen thus: given \eqref{H-cnd}, Proposition~\ref{supp-prop} 
implies that any
$\xi+\eta$ in $\supp\cal F(a_{k-j,k-l}(x,D)u_{k-l})$ for 
large $k$ fulfils 
\begin{equation}
  |\xi+\eta|\ge \tfrac{1}{C} |\eta|-1\ge \tfrac{11}{20C}2^{k-l}-1
  \ge (\tfrac{11}{40C}-2^{-k})2^k>\tfrac{1}{4C}2^k.
\end{equation}
To complete the proof of Theorem~\ref{cont-thm} one can modify the estimates
\eqref{a1-ineq} ff.\ into $L_p(\ell^s_q)$ estimates; then
convergence criteria for series of distributions, eg Theorems~3.6--3.7 of
\cite{Y1}, apply by \eqref{S1}--\eqref{S2} (like arguments
used in \cite{Mey80,Y1,Mar96} etc.). The ball on the r.h.s.\ of \eqref{S2}
only yields estimates of $\norm{a^{(2)}(x,D)u}{F^s_{p,q}}$ for $s>0$, as
is well known. But if \eqref{H-cnd} holds, one can, by \eqref{S2'}, use the
criteria for series with spectra in dyadic annuli, like for $a^{(1)}$ and
$a^{(3)}$ (the finitely many other terms of $a^{(2)}$ are in
$\bigcap_{s>0}F^s_{p,q}$). 

\begin{remark}
The class
$\op{OP}(S^d_{1,1}(\Rn\times\Rn))$ was first treated in 
$F^{s}_{p,q}$-spaces by Runst \cite{Run85ex}, but unfortunately
the proofs are somewhat flawed, since in 
Lemma~1 there the spectral estimates require a support rule under rather
weak assumptions, like in Prop.~\ref{supp-prop} above.
This was seemingly overlooked in
\cite{Run85ex} and by Marschall \cite{Mar96}.
Using the $\varphi$-decomposition of Frazier and Jawerth~\cite{FJ2},
Torres \cite{Tor90} 
extended the $H^s_p$-continuity of \cite{Mey80} to the
$F^s_{p,q}$-scale. 
The borderline $s=0$ was treated by Bourdaud
\cite[Thm.~1]{Bou88}; his result on $B^0_{p,1}$ is improved by
Thm.~\ref{Fp1-thm} above. 
Thm.~\ref{cont-thm} is a novelty concerning \eqref{H-cnd}. 
\end{remark}

\end{document}